\documentclass[conference]{IEEEtran}

\IEEEoverridecommandlockouts

\usepackage{cite}
\usepackage{amsmath,amssymb,amsfonts}
\usepackage{algorithmic}
\usepackage{graphicx}
\usepackage[export]{adjustbox}
\usepackage{booktabs}
\usepackage{textcomp}
\usepackage{xcolor}
\def\BibTeX{{\rm B\kern-.05em{\sc i\kern-.025em b}\kern-.08em
    T\kern-.1667em\lower.7ex\hbox{E}\kern-.125emX}}
\usepackage{url}
\usepackage[hidelinks]{hyperref}

\begin{document}

\title{\LARGE Optimal Scheduling of Multi-Site AI Data Centers in Day-Ahead Markets}

\author{
  \IEEEauthorblockN{Xinyang Zhou and Guido Cavraro}
  \IEEEauthorblockA{National Laboratory of the Rockies, Golden, CO, USA; Email: xinyang.zhou@nlr.gov, guido.cavraro@nlr.gov}
}
\maketitle

\begin{abstract}
AI data centers are among the fastest-growing loads on the power grid, and their flexibility is increasingly sought by grid operators. Using it requires a model that links what an operator controls, namely GPU power caps, job timing and job placement, to the power drawn at each site and the electricity prices it faces. This paper develops such a framework for a multi-site AI data-center fleet. It represents GPU power-cap levels; deferrable jobs with deadlines and inference with response-time headroom; and site facility power, and casts the fleet schedule as a linear program solved on day-ahead prices. A year-long simulation on 2025 California Independent System Operator (CAISO) day-ahead prices in three data-center markets demonstrates the framework. Coordinated scheduling lowers the fleet's energy cost by 6.3\% at equal work. Four operating modes, each adding one capability, split this saving into efficient GPU operation that needs no price signal, price response in time, and routing work between sites.
\end{abstract}


\section{Introduction}\label{sec:intro}
U.S.\ data centers used 4.7\% of the nation's electricity in 2024 and may use 9.5--15.3\% by 2030~\cite{lbnl2026}, driven mainly by AI servers. About 100~GW of new large load could connect if curtailed by about 0.5\% of its maximum annual energy use~\cite{norris2025}. Field tests have held a GPU cluster 25\% below baseline for three hours during a grid peak~\cite{colangelo2025} and, in a later London trial, cut a cluster's demand by up to 40\% and followed reduction requests lasting up to 10~h~\cite{nationalgrid2026}. Much of an AI fleet's work can be paused, capped or deferred~\cite{radovanovic2023}, and an operator with several sites can also move it~\cite{liu2011greening}.

Turning this potential into grid value requires a model that connects the operator's controls to the grid. On the data-center side, the controls act on individual GPUs and jobs: a power cap trades throughput for power, training can pause and resume from checkpoints, offline jobs have deadlines, and inference must meet response-time targets. On the grid side, what matters is each site's facility power, including cooling, and the locational prices and market rules it faces. Geographic load balancing~\cite{liu2011greening} and carbon-aware scheduling~\cite{radovanovic2023} move work between sites and over time, respectively, but do not model GPU power caps. Grid studies, conversely, typically assume that a fixed share of data-center load can be deferred~\cite{mazumder2026,riepin2025}; a trace of 155,410 GPUs instead derived flexibility from the workloads eligible to curtail and found that it averages about 6\% of facility power \cite{li2026beyond}.

Such a model must also be clear about where savings come from. An AI fleet that faces electricity prices can lower its bill in three ways. It can run GPUs at lower power caps, which finishes the same work with less energy at any price. It can shift work in time, toward hours when power is cheap. And it can move work between sites, toward whichever site has cheaper power at that moment. Only the last two respond to prices, yet all three appear as savings when an optimized schedule is compared with a fleet that runs every job at full power on arrival, which counts efficiency as flexibility.

This paper makes the following contributions:
\begin{enumerate}
  \item A scheduling framework for a multi-site AI data-center fleet that links GPU power-cap levels and workload constraints to site facility power and locational electricity prices. It is a linear program that sets each site's hourly power a day ahead from day-ahead (DA) prices.
  \item A year-long demonstration on 2025 CAISO DA prices in three data-center markets. Four operating modes, each adding one capability, split the saving into efficiency, price response in time and routing between sites, with confidence intervals. Coordinated scheduling saves 6.3\% at equal work; efficiency and price response carry nearly all of it, and routing adds a small share, mainly on days when prices at the sites diverge.
\end{enumerate}
Because the framework ties GPU-level controls to grid-facing quantities, it provides a basis for valuing data-center grid services such as demand response and reserves.

The rest of this paper is organized as follows. Section~\ref{sec:model} models the fleet and its operating modes. Section~\ref{sec:setup} describes the simulation setup, Section~\ref{sec:results} presents the results, and Section~\ref{sec:conclusion} concludes this work.

\section{Fleet Scheduling Model}\label{sec:model}

A data center fleet operator manages a fleet of data centers over time windows of length $\Delta$ indexed by $t$, so that the $t$-th time window corresponds to the time interval $[(t-1)\Delta, t\Delta]$. Each data center $i$ hosts $N_i$ identical GPUs. The fleet operator sets for every time interval 1) how many GPUs are assigned to each computational job, 2) at which power cap the GPUs run, and, 3) when sites are coordinated, to which site movable work is transferred.

We assume for simplicity that the fleet is a \emph{price taker}, i.e., the fleet operator schedules against locational marginal prices (LMPs) which are not affected by the data center loads.
We will cast the fleet operator scheduling program as a linear program (LP), which can be efficiently solved by off the shelf methods, with continuous decision variables.
The formulation follows geographic load balancing~\cite{liu2011greening}, with per-job deadlines in place of the daily completion window of~\cite{radovanovic2023}.

We assume two types of work here.
\emph{Deferrable jobs}, indexed by $m$, must finish by a deadline but may be paused, slowed and, in some cases, moved to other sites. They range from training runs that last weeks, use thousands of GPUs and are pinned to their site, to short jobs such as fine-tuning, evaluation and data processing that last hours and may run at any site. \emph{Inference} serves live user requests that cannot wait; it is indexed as $m=\mathrm{inf}$ so that sums over $m$ in the following include it.

\begin{figure*}[t]
  \centering
  \includegraphics[trim={0 0 0 2mm}, clip, width=0.8\textwidth]{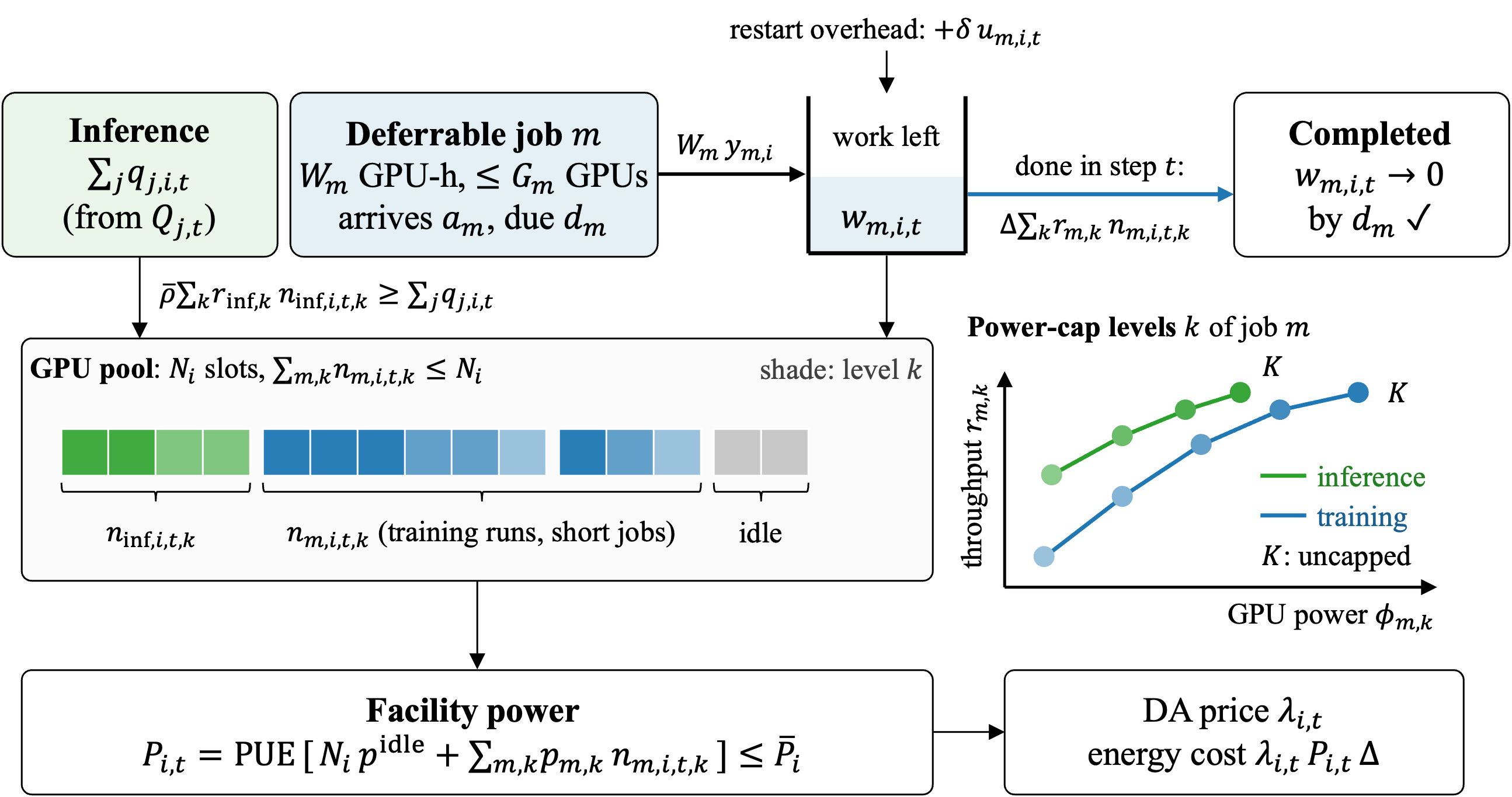}
  \caption{The scheduling model at one site $i$. Inset: throughput against GPU power at the levels of Table~\ref{tab:param}.}
  \vspace{-6pt}
  \label{fig:model}
\end{figure*}

In the following, Sections~\ref{sec:levels}--\ref{sec:obj} build the LP; Section~\ref{sec:control} describes the day ahead (DA) schedule and the four operating modes we compare. Fig.~\ref{fig:model} shows how the pieces fit together at one site, and Table~\ref{tab:nom} summarizes the notation.

\begin{table}[t]
\centering
\caption{Notation}\label{tab:nom}
\begin{tabular}{@{}l p{0.74\columnwidth}@{}}
\toprule
\multicolumn{2}{@{}l}{\textbf{Indices}}\\
$i$, $j$; $t$, $\Delta$, $H$ & sites; time step, its length, horizon end\\
$k$, $K$ & power-cap level; uncapped (natural) level\\
$m$ & deferrable job (training run or short job), or inference\\
\multicolumn{2}{@{}l}{\textbf{Variables}}\\
$n_{m,i,t,k}$ & GPUs of job $m$ at site $i$ in step $t$ at level $k$\\
$u_{m,i,t}$ & GPUs of job $m$ that start up in step $t$\\
$y_{m,i}$ & share of job $m$ placed at site $i$\\
$q_{j,i,t}$ & inference demand from site $j$ served at site $i$\\
$w_{m,i,t}$; $P_{i,t}$ & remaining work (GPU-h); facility power (MW)\\
\multicolumn{2}{@{}l}{\textbf{Parameters}}\\
$\phi_{m,k}$, $r_{m,k}$ & GPU power (TDP share), relative throughput of $m$ at $k$\\
$p^{\mathrm{idle}}$, $p_{m,k}$ & power of an idle GPU slot; extra power of a busy one\\
$p^{\mathrm{TDP}}$, $p^{0}$ & GPU TDP; fixed extra power of a busy slot\\
$N_i$, $\bar P_i$ & number of GPUs and power limit of site $i$\\
$W_m$, $G_m$, $s_m$ & work at arrival (GPU-h), GPU limit, cost to move (\$)\\
$i_m$, $a_m$, $d_m$ & home site, arrival step, deadline of job $m$\\
$\delta$ & work lost per restarted GPU (h)\\
$Q_{j,t}$; $\bar\rho$ & inference demand arising at site $j$; utilization cap\\
$\mathrm{PUE}$ & facility power divided by IT power\\
$\lambda_{i,t}$ & DA electricity price (\$/MWh)\\
$\pi^{\mathrm{late}}$, $\pi^{\mathrm{route}}$ & prices of lateness and inference routing\\
$v_m$; $\kappa$ & value of leftover work; value of GPU time\\
\bottomrule
\end{tabular}
\end{table}

\subsection{GPU Power Levels}\label{sec:levels}
Power capping~\cite{patel2024polca} lets a GPU run at a few discrete levels $k$,
which depend on its workload class. A GPU assigned to job $m$ at level $k$ draws a fraction $\phi_{m,k}$ of its thermal design power (TDP) and has corresponding relative throughput $r_{m,k}$: one hour at level $k$ completes $r_{m,k}$ GPU-hours of uncapped work. The last level $K$ is the natural, uncapped draw, with $r_{m,K}=1$.
We assume a GPU can also split a step between two levels; it then averages their power and throughput and operates on the straight line between them; see the inset of Fig.~\ref{fig:model}. Throughput is therefore a piecewise-linear function of power. The levels show diminishing returns: each extra watt buys less throughput, so this function is concave. This allows us to trade off GPU throughput for power saving, especially in presence of flexible job deadlines.
We define $n_{m,i,t,k}$ as the number of GPUs allocated to job $m$ at site $i$ and time interval $t$ capped at the power level $k$. For simplicity, we consider $n_{m,i,t,k}$ as a continuous rather than discrete variable, given the big number of GPUs (tens of thousands) in each site.

\subsection{Deferrable Jobs}\label{sec:work}
Deferrable jobs are measured in GPU-hours (GPU-h) at uncapped throughput: $W=1{,}000$~GPU-h keeps 1{,}000 uncapped GPUs busy for one hour, or 1{,}000 GPUs at $r_{m,k}=0.5$ for two hours. A deferrable job $m$ that arrives is characterized by its work $W_m$ and by the maximum number of GPUs that may be allocated to it at once  $G_m$.
To each site $i$, the schedule assigns a share $y_{m,i}$ of job $m$, with $\sum_i y_{m,i}=1$, and allocates $n_{m,i,t,k}$ GPUs at level $k$. Usually, training runs are pinned to their home site $i_m$, short jobs are more easily moved.
After job $m$ started to be processed in site $i$, the remaining work left at time $t$ is denoted as $w_{m,i,t}$ and evolves like an account balance (the tank in Fig.~\ref{fig:model}) over the planning horizon, which ends at step $H$:
\begin{subequations}\label{eq:work}
{\medmuskip=2mu\thickmuskip=3mu
\begin{align}
w_{m,i,t+1} &= w_{m,i,t}-\Delta \textstyle\sum_k r_{m,k} n_{m,i,t,k}+\delta u_{m,i,t}, \label{eq:work-bal}\\
u_{m,i,t} &= \max\{\textstyle\sum_k(n_{m,i,t,k}-n_{m,i,t-1,k}),~0\}, \label{eq:work-up}\\
w_{m,i,a_m} &= W_m\, y_{m,i}, \label{eq:work-init}\\
\textstyle\sum_k n_{m,i,t,k} &\le G_m\, y_{m,i}, \label{eq:work-gpu}\\
w_{m,i,H} &\le W_m\, y_{m,i}\,\max\{d_m-H, 0\}/(d_m-a_m), \label{eq:work-pace}
\end{align}}
\end{subequations}
where $u_{m,i,t}$ is the GPUs of job $m$ that start up in step $t$, and $a_m$ and $d_m$ are the arrival step and deadline of job $m$.
Equation~\eqref{eq:work-bal} updates the remaining work by subtracting the work GPUs complete and by adding the cost of restarting GPUs: a GPU that resumes a paused job must reload the last checkpoint and redo the work done after last saving, so it spends $\delta$ hours busy without progress.
Equation~\eqref{eq:work-up} counts the GPUs that start up: the increase in running GPUs. The LP imposes~\eqref{eq:work-up} linearly as $u_{m,i,t}\ge\sum_k(n_{m,i,t,k}-n_{m,i,t-1,k})$ and $u_{m,i,t}\ge 0$: since restarts only add work, $u_{m,i,t}$ meets this bound at the optimum.
Constraints~\eqref{eq:work-init} and~\eqref{eq:work-gpu} give site $i$ its share of the job's work and of its GPU limit.
Constraint~\eqref{eq:work-pace} is a pace limit:
the work left at the end of the horizon may not exceed a straight-line path to the deadline. A training run due in 10 days must finish at least 10\% of its work on the first day, and a short job due within the horizon must finish completely.

\subsection{Online Inference}\label{sec:inf}
Inference is live traffic and cannot be delayed, so it has no work balance. The demand arising at site $j$, $Q_{j,t}$, is measured in GPU-equivalents, the number of uncapped GPUs it would fully occupy during the time slot $t$. A part $q_{j,i,t}$ of it is served at site $i$ by $n_{\mathrm{inf},i,t,k}$ inference GPUs at level $k$. At this DA planning stage, we assume the serving GPUs are loaded to at most a fraction $\bar\rho$ of their capacity, the same at every site, which leaves slack for bursts of requests~\cite{liu2011greening,stojkovic2025dynamollm}. The inference constraints are:
\begin{subequations}\label{eq:inf}
\begin{align}
\textstyle\sum_i q_{j,i,t} &= Q_{j,t}, \label{eq:inf-dem}\\
\textstyle\bar\rho\sum_k r_{\mathrm{inf},k} n_{\mathrm{inf},i,t,k} &\ge \textstyle\sum_j q_{j,i,t}, \label{eq:inf-cap}\\
\textstyle\sum_{m,k}n_{m,i,t,k} &\le N_i. \label{eq:inf-site}
\end{align}
\end{subequations}
Constraint~\eqref{eq:inf-dem} serves all inference demand, \eqref{eq:inf-cap} gives each site enough serving GPUs for the demand it takes at utilization $\bar\rho$, and \eqref{eq:inf-site} keeps all GPUs in use (including deferrable jobs and inference) within the $N_i$ GPUs of the site. In fleet mode demand may be served at any site; otherwise $q_{j,i,t}=0$ for $i\neq j$. With $\bar\rho=0.8$, serving 1{,}000 GPU-equivalents takes 1{,}250 uncapped GPUs, or 1{,}543 at 46\% of TDP. The lower cap may need more GPUs but less power overall.

\subsection{Facility Power}\label{sec:power}
Every GPU slot, with its share of the server, draws $p^{\mathrm{idle}}$ even when idle; a busy slot of job $m$ at level $k$ draws
\begin{equation}
p_{m,k}=p^{\mathrm{TDP}}\phi_{m,k}+p^{0} \label{eq:pmk}
\end{equation}
more, where $p^{\mathrm{TDP}}\phi_{m,k}$ is the GPU's power at level $k$ and $p^{0}$ is the server's extra power when busy, net of the GPU's idle draw.
Facility power is IT power times the power usage effectiveness (PUE), which adds cooling and power-conversion losses:
\begin{equation}
P_{i,t}=\mathrm{PUE}\Big[N_i p^{\mathrm{idle}}+\textstyle\sum_{m,k}p_{m,k} n_{m,i,t,k}\Big]\le \bar P_i. \label{eq:power}
\end{equation}
The first term in brackets is fixed; the second depends on the schedule. We take the PUE as a constant, $\mathrm{PUE}=1.1$, near the fleet averages reported by hyperscale operators~\cite{google_pue}.

\subsection{Objective and Optimization Problem}\label{sec:obj}
The goal for the operator is to minimize the following objective function:
\begin{align}
J=\;&\textstyle\sum_{i,t}\lambda_{i,t}P_{i,t}\Delta
+\pi^{\mathrm{late}}\sum_{m,i,t\ge d_m}w_{m,i,t}\Delta \nonumber\\
&\textstyle+\sum_{m,i}v_m\,w_{m,i,H}
+\sum_{m}\sum_{i\neq i_m}s_m\, y_{m,i} \nonumber\\
&\textstyle+\pi^{\mathrm{route}}\sum_{t}\sum_{j,\,i\neq j}q_{j,i,t}\Delta, \label{eq:obj}
\end{align}
where $\lambda_{i,t}$, $\pi^{\mathrm{late}}$, and $\pi^{\mathrm{route}}$ are price for electricity, lateness, and inference routing.
The first term captures the energy bill, the second term models a late penalty $\pi^{\mathrm{late}}$ on each GPU-hour of work still unfinished after its deadline, the third term sums the value of the work left at the end of the horizon, which must still be done later, the fourth term is the cost of moving a job's data away from its home site, and the fifth term is a small routing cost for serving inference away from where it arises. Here, the leftover work is valued at:
\begin{equation}
v_m=p_{m,K}\,\mathrm{PUE}\,\bar\lambda+\kappa, \label{eq:value}
\end{equation}
where the first term is the energy needed to finish the work later, at this day's mean price $\bar\lambda$, and the second, $\kappa$ in \$/GPU-h, is the value of the GPU time the work will occupy, since a GPU-hour deferred today displaces other work later. The two differ greatly in scale: at \$39/MWh, the energy of a training GPU-hour costs about \$0.03, while an H100 GPU rents for about \$2/h.
\$2/h is near the rental value, so training is never deferred, consistent with~\cite{li2026beyond}, and every mode completes the same work.

We can now cast our scheduling problem as the LP:
\begin{subequations}
\begin{align}
\min_{n,\,u,\,w,\,y,\,q,\,P\,\ge\,0} &\quad J \\
\text{s.t. } & \quad \textstyle\sum_i y_{m,i}=1, \label{eq:lp} \\
& \quad \eqref{eq:work}\text{--}\eqref{eq:power},\;\;
\end{align}
\end{subequations}
with $y_{m,i_m}=1$ for training runs and in the site-level modes. The LP input parameters are the prices $\lambda$, the workload $(W_m,Q_{j,t})$ and the starting state; its outputs are the GPU counts per level, the placement and the facility power at each site. To keep the LP always feasible, the implementation softens the pace limit~\eqref{eq:work-pace} and the inference balance~\eqref{eq:inf-dem} with penalties of \$20/GPU-h and \$5{,}000 per GPU-equivalent-hour; neither is incurred in the simulation results.

\subsection{DA Schedule and Operating Modes}\label{sec:control}
The fleet schedules each operating day in CAISO's DA market.
The LP is solved once, in hourly steps ($\Delta=1$h) over the 24 hours of the day, with the DA LMPs $\lambda_{i,t}$ and the day's workload. Its hourly power $P_{i,t}$ is the fleet's DA purchase, which it then consumes as scheduled, so its bill is the first term of~\eqref{eq:obj}, $\sum_{i,t}\lambda_{i,t}P_{i,t}\Delta$.

We compare four modes. Each adds one capability to the one before it, so that the difference between consecutive modes isolates one source of saving:
\begin{enumerate}
  \item \emph{Uncoordinated}: every job runs on arrival at its home site at full power, inference first, then training runs, then short jobs in arrival order.
  \item \emph{Price-blind}: the site-level model is solved at a flat price, with a negligible preference for early work. This isolates operational efficiency.
  \item \emph{Site-level}: the model is solved on LMPs, with all work and inference kept at its home site.
  \item \emph{Fleet-coordinated}: the model is solved on LMPs and also optimizes routing ($y$ and $q$).
\end{enumerate}
All modes pay the DA LMPs for the power they draw. Thus (b) minus (a) measures efficiency, (c) minus (b) measures price response in time, and (d) minus (c) measures routing.

\section{Simulation Setup}\label{sec:setup}
\begin{table}[t]
\centering
\caption{Main parameters}
\label{tab:param}
\setlength{\tabcolsep}{3pt}
\begin{tabular}{@{}p{0.47\columnwidth}p{0.50\columnwidth}@{}}
\toprule
Parameter & Value\\
\midrule
GPUs per site & 90k Silicon Valley; 60k LA, San Diego\\
GPU TDP $p^{\mathrm{TDP}}$; GPU idle; server busy/idle & 700; 70; 450/250~W\\
$p^{\mathrm{idle}}$; $p^{0}$ (derived) & 320; 130~W\\
Training power-cap levels $(\phi,r)$ & $(0.85,1)$, $(0.75,0.96)$, $(0.65,0.88)$, $(0.55,0.76)$, $(0.45,0.62)$\\
Inference power-cap levels $(\phi,r)$ & $(0.70,1)$, $(0.63,0.96)$, $(0.55,0.90)$, $(0.46,0.81)$ \cite{apartsin2026joule}\\
Short-job power-cap levels $(\phi,r)$ & $(0.80,1)$ \ldots $(0.40,0.70)$\\
Inference utilization cap $\bar\rho$; PUE & 0.8; 1.1 \cite{google_pue}\\
Restart overhead $\delta$, training/short & 0.25/0.05~hour\\
$\pi^{\mathrm{late}}$; $s_m$ & \$10/GPU-h per h; \$0.01/GPU-h moved\\
$\pi^{\mathrm{route}}$; $\kappa$ & \$0.05/GPU-h; \$2/GPU-h\\
\bottomrule
\end{tabular}
\end{table}

\subsection{Fleet and Prices}
Three sites host 210,000 GPUs in total. When run uncoordinated, they average about 190~MW with a mean daily peak of 230~MW.
The sites are placed in CAISO's three largest data-center markets, Silicon Valley, the Los Angeles basin and San Diego. Each site pays the LMP that CAISO publishes for its local area (PG\&E South Bay, SCE Core and SDG\&E).

The 2025 hourly DA LMPs come from CAISO's OASIS data portal~\cite{caiso_oasis}. DA prices at the three sites averaged \$39/MWh and peaked at \$199/MWh, and midday solar output drove them below zero at some site on 102 days. Silicon Valley was the most expensive market, with an average price of \$43/MWh against \$35/MWh in Los Angeles.

\subsection{Workload and Runs}
We assume each site runs two training runs on 26\% of its GPUs, with 5--30 days left to run and a 15\% deadline slack. These runs are pinned to their cluster, as large synchronous runs are. Short jobs arrive following a business-day profile, totaling 19\% of GPU capacity per day with lognormal volume, and may move between sites; those arriving in the same 3-h window with the same deadline are modeled as one job. Each short job is due 6, 12 or 24~h after arrival, with assumed probabilities 0.2, 0.3 and 0.5, within the 24-h window. Online inference follows a diurnal shape with an assumed peak-to-trough ratio of 1.6 for the fleet's mix of services.
Average GPU allocation is about 75\%.

The workload is taken as known a day ahead. Each day is simulated from a fresh state with an independent workload draw, and training runs enter at a random point in their run.
Every day of 2025 (365 days) is simulated in all four modes.
Models are built in Pyomo and solved with HiGHS.

\section{Simulation Results}\label{sec:results}
\subsection{Efficiency, Price Response and Routing}\label{sec:decomp}
\begin{table}[t]
\centering
\caption{Totals and averages of energy and cost in 2025.}
\label{tab:year}
\begin{tabular}{@{}l rrrr@{}}
\toprule
 & Uncoord. & Price-blind & Site & Fleet\\
\midrule
Energy cost (\$M) & 65.32 & 62.87 & 61.36 & 61.24\\
Saving vs.\ uncoord.\ (\%) & -- & 3.76 & 6.07 & 6.25\\
Energy (GWh) & 1{,}676.1 & 1{,}597.6 & 1{,}617.6 & 1{,}617.6\\
Avg.\ price (\$/MWh) & 38.97 & 39.35 & 37.93 & 37.86\\
Mean daily peak (MW) & 229.2 & 192.0 & 216.2 & 216.4\\
\bottomrule
\end{tabular}
\end{table}

\begin{figure*}[t]
  \centering
  \includegraphics[width=1\textwidth]{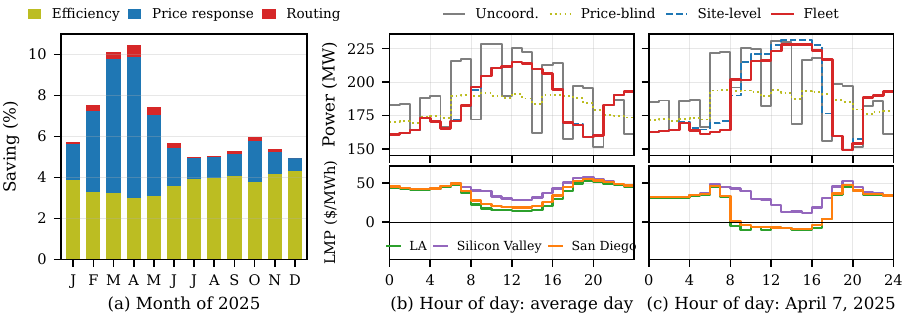}
  \caption{Results over 2025. (a) Monthly fleet saving relative to the uncoordinated mode.
  (b) Average day: fleet power by mode and DA LMP by site, averaged over all days for each hour; the site-level and fleet curves nearly coincide. (c) April 7:
  the price-responsive modes fill the negative-price solar hours, and the fleet mode also moves load out of Silicon Valley.}
  \vspace{-6pt}
  \label{fig:results}
\end{figure*}

Table~\ref{tab:year} summarizes the year. All modes complete the same 835 million GPU-hours of deferrable work, serve all inference demand and finish every short job by its deadline, so the comparison is at equal output. Fleet coordination cuts energy cost by 6.3\%,
or about \$4.1M per year and \$2.5/MWh consumed. Taking differences between consecutive modes splits this saving into three parts, each in percentage points of the uncoordinated cost; Fig.~\ref{fig:results}(a) shows them by month.

\emph{Efficiency: 3.8\%.} Run at a flat price, the optimizer spreads work across spare GPUs at lower power caps, cutting energy by 4.7\% and the mean daily peak by 16\%. The gain comes from inference and short jobs. At 46\% of TDP, a GPU runs each inference request about 1.2 times slower but draws much less power.
For training, the gain is under 1\%: its runs already use their full GPU limit, so a lower cap would only slow them.

\emph{Price response: 2.3\%.} Responding to LMPs lowers the average price paid by 3.6\% relative to the price-blind schedule, mainly by moving work into the midday hours, when solar output lowers prices (Fig.~\ref{fig:results}(b)).
Its value follows the solar season: about 7\% in March and April, when midday prices fell below zero on most days, but only about 1\% from July to September.

\emph{Routing: 0.2\%.}
It exceeds 1\% on 23 days. On April 7 (Fig.~\ref{fig:results}(c)), for example, midday DA prices fell below zero in Los Angeles and San Diego but averaged \$25/MWh in Silicon Valley. The fleet mode moved about 6~MW of midday load out of Silicon Valley, which was worth 2.1\% that day.

\section{Conclusion and Future Work}\label{sec:conclusion}
This paper developed a modelling and simulation framework that links an AI data-center fleet's GPU power caps, job timing and job placement to site power and locational prices. Formulated as a linear program, it remains tractable at fleet scale, and nested scenarios that release one degree of freedom at a time separate efficiency from price-responsive flexibility. On 2025 CAISO day-ahead prices, coordinated scheduling lowers energy cost by 6.3\% at equal work: 3.8\% from price-independent power-cap efficiency, 2.3\% from temporal shifting and 0.2\% from spatial routing. Routing gains little within one ISO, where nodal prices are strongly correlated, but should be worth more across sites with larger price differences, such as congested nodes or different ISOs.

Ongoing work will add real-time operation, re-solving the LP every 15~min by model-predictive control on CAISO's fifteen-minute market prices, accounting for prediction errors and uncertainties. Future work will use the model to schedule and value grid services such as demand response, voltage and frequency regulation, and couple it with the network through congestion-aware allocation and AC optimal power flow.

\section*{Acknowledgment}
This work was authored in part by the National Laboratory of the Rockies for the U.S. Department of Energy (DOE) under Contract No. DE-AC36-08GO28308. Funding provided by the Advanced Modeling Grid Research Program at the Office of Electricity, U.S. Department of Energy. The views expressed in the article do not necessarily represent the views of the DOE or the U.S. Government. The U.S. Government retains and the publisher, by accepting the article for publication, acknowledges that the U.S. Government retains a nonexclusive, paid-up, irrevocable, worldwide license to publish or reproduce the published form of this work, or allow others to do so, for U.S. Government purposes.

The authors used Anthropic’s Claude Opus 5.5 to revise
some text in the paper. All AI-assisted content was reviewed,
edited, and verified by the authors, who take full responsibility for the paper.

\bibliographystyle{IEEEtran}
\bibliography{refs}

\end{document}